\documentclass[10pt]{article}
\usepackage{amsthm,amsmath,amsfonts}
\usepackage[authoryear]{natbib}
\usepackage{color,graphicx,rotating}
\usepackage{setspace}
\usepackage{bm}
\usepackage[sf,bf,tiny]{titlesec}
\titlelabel{\thetitle. \enspace}
\date{}
\title{\normalsize\bf  The $\bm{\Lambda}$-Fleming-Viot process and a connection with Wright-Fisher diffusion.}
\author{{\sc Robert C. Griffiths
\thanks{\noindent
Professor R. C. Griffiths, Department of Statistics, University of Oxford, 1 South Parks Rd, Oxford, OX1 3TG, UK; Ph +44 1865 281237;Fax +44 1865 272595;
email  griff@stats.ox.ac.uk; \emph{Version 7c}}
}\\\medskip
\emph{Oxford University} }
\begin{document}
\maketitle
{\small
\section*{Abstract}
The $d$-dimensional $\Lambda$-Fleming-Viot generator acting on functions $g(\bm{x})$, with $\bm{x}$ being a vector of $d$ allele frequencies, can be written as a Wright-Fisher generator acting on functions $g$ with a modified random linear argument of $\bm{x}$ induced by partitioning occurring in the $\Lambda$-Fleming-Viot process. The eigenvalues and right polynomial eigenvectors are easy to see from this representation. The two-dimensional process, which has a one-dimensional generator, is considered in detail. A non-linear equation is found for the Green's function.
In a model with genic selection a proof is given that there is a critical selection value such that if the selection coefficient is greater or equal to the critical value then fixation, when the boundary 1 is hit, has probability 1 beginning from any non-zero frequency. This is an analytic proof different from proofs by \cite{DEP2011} and \citet{F2013}. 
An application in the infinitely-many-alleles $\Lambda$-Fleming-Viot process is finding an interesting identity for the frequency spectrum of alleles that is based on size-biassing.
The moment dual process in the Fleming-Viot process is the usual $\Lambda$-coalescent tree back in time. The Wright-Fisher representation using a different set of polynomials $g_n(x)$ as test functions produces a dual death process which has a similarity to the Kingman coalescent and decreases by units of one. The eigenvalues of the process are analogous to the Jacobi polynomials when expressed in terms of $g_n(x)$, playing the role of $x^n$.
$\mathbb{E}\big [g_n(X)\big ]$ under the stationary distribution when there is mutation is analogous to the $n^{\text{th}}$ moment in a Beta distribution. There is a $d$-dimensional version $g_{\bm{n}}(\bm{X})$, and even an intriguing Ewens' sampling formula analogy when $d \to \infty$.
\medskip

\noindent
\emph{Keywords}: $\Lambda$-coalescent Fleming-Viot process, Wright-Fisher diffusion processs.

\noindent
2010 Mathematics Subject Classification: Primary 60G99, Secondary 92D15
}
%%%%%
\newpage
\section{Introduction}
The $d$-dimensional $\Lambda$-Fleming-Viot process $\{\bm{X}_t\}_{t\geq 0}$ representing frequencies of $d$ types of individuals in a population has state space
$
\Delta = \{\bm{x} \in [0,1]^d: \sum_{i\in [d]}x_i \leq 1\}
$
with generator ${\cal L}$ acting  on functions in $C^2(\Delta)$ described by
\begin{equation}
{\cal L}g(\bm{x}) = \int_0^1\sum_{i=1}^dx_i\big (g(\bm{x}(1-y)+y\bm{e}_i) - g(\bm{x})\big )\frac{\Lambda(dy)}{y^2}.
\label{L:basic}
\end{equation}
In general $\Lambda$ is a non-negative finite measure on $[0,1]$. We take a time scale so that $\Lambda \equiv F$ is a probability measure on $[0,1]$. Informally the population is partitioned at events of change by choosing type $i \in \{1,2,\ldots,d\}$ to reproduce with probability $x_i$, then rescaling the population with additional offspring $y$ of type $i$ so that the frequencies are $\bm{x}(1-y) + y\bm{e}_i$, at rate $y^{-2}F(dy)$. 
If $F$ has a single atom at $0$, then $\{\bm{X}_t\}_{t\geq 0}$ is the $d$-dimensional Wright-Fisher diffusion process in $\Delta$ with generator
\begin{equation}
\label{L:WF}
{\cal L} = \frac{1}{2}\sum_{i,j=1}^dx_i(\delta_{ij}-x_j)\frac{\partial^2}{\partial x_i\partial x_j}
\end{equation}
acting on functions in $C^2(\Delta)$.
The general process $\{\bm{X}_t\}_{t \geq 0}$ with generator (\ref{L:basic}) has a Wright-Fisher diffusive component if $F({0})>0$ and discontinuous sample paths from jumps where the frequencies are changed by adding mass $y$ from the points of $F$ in $(0,1]$ to the frequency of a type and rescaling the frequencies. Eventually the process becomes absorbed into one state in $\{\bm{e}_i\}_{i=1}^d$.
\citet{EW2006} introduced a model where $F$ has a single point of increase in $(0,1]$ with a possible atom at zero as well. A natural class that arises from discrete models are Beta-coalescents, particularly when $F$ has a Beta$(2-\alpha,\alpha)$ density coming from a discrete model where the offspring distribution tails are asymptotic to a power law of index $\alpha$. This Beta-coalescent model is studied in \citet{S2003,BBCEMSW2005}. \citet{BB2009} describe the $\Lambda$-Fleming-Viot process and discrete models whose limit gives rise to the process.

The $\Lambda$-coalescent is a random tree back in time which has multiple merger rates for a specific $2 \leq k\leq n$ edges merging while $n$ edges in the tree of
\begin{equation}
\lambda_{nk} = \int_0^1x^k(1-x)^{n-k}\frac{\Lambda (dx)}{x^2},\>k \geq 2.
\label{lnk:rates}
\end{equation}
After coalescence there are $n-k+1$ edges in the tree.
The process is often regarded as having a state space on the set of partitions $\Pi_\infty$ of the positive integers. The leaves of an infinite leaf $\Lambda$-coalescent tree at time $t=0$ are labelled with singleton sets 
$\{1\},\{2\},\ldots$ and edges at time $t$ are labelled by sets in $\Pi_\infty(t)$. The number of blocks at time $t$ is the number of sets in the partition $\Pi_\infty(t)$, denoted by $|\Pi_\infty(t)|$, which is the same as the number of edges in the tree at time $t$. If there are $n$ edges at time $t$, and $k$ merge at $t^{+0}$, then a new partition is formed by taking the 
union of the $k$ partition blocks in the merger for the parent block at $t^{+0}$. This occurs at rate $\lambda_{nk}$. The $\Lambda$-coalescent is said to \emph{come down from infinity} if for all $t>0$, $P(|\Pi_\infty(t)| < \infty) = 1$, which is equivalent to an infinite-leaf $\Lambda$-coalescent tree at $t=0$ having a finite number of edges at any time $t>0$ back with probability 1.

The $\Lambda$-coalescent process was introduced by \citet{DK1999,P1999,S1999} and has been extensively studied \citep{P2002,B2009}. The coalescent process is a moment dual to the $\Lambda$-Fleming-Viot process. See for example \citet{E2012}. There is a distinction between an untyped coalescent process and a typed process such as in \citet{EGT2010}. 

There is a connection between continuous state branching processes and the $\Lambda$-coalescent. For example see \citet{BLG2003,BLG2006,BBCEMSW2005,BBLa2012,BBLb2012}. The connection is through the Laplace exponent
\begin{equation}
\psi(q) = \int_0^1 \big (e^{-qy} -1 +qy\big )y^{-2}\Lambda(dy).
\label{Laplaceexponent}
\end{equation}
\citet{BLG2006} showed that the $\Lambda$-coalescent comes down from infinity under the same condition that the continuous state branching process becomes extinct in finite time, that is when
\begin{equation}
\int_1^\infty\frac{dq}{\psi(q)} < \infty.
\label{comedown:0}
\end{equation}
\citet{S2000} proved earlier that coming down from infinity was equivalent to
\begin{equation*}
\sum_{n=2}^\infty \Biggl [\sum_{k=2}^n(k-1){n\choose k}\lambda_{nk}\Biggr ]^{-1} < \infty.
%\label{comedown:a}
\end{equation*}

In this paper we express the $\Lambda$-Fleming-Viot generator acting on functions as a Wright-Fisher diffusion generator where the argument of the function is replaced by a random linear transformation. For example if $d=2$ the generator acting on functions of $x_1=x$ in ${\cal C}^2([0,1])$ is specified by
\begin{equation}
{\cal L}g(x) = \int_0^1\Big [x\big (g(x(1-y)+y\big )-g(x)\big )+(1-x)\big (g(x(1-y)) - g(x)\big )\Big ]\frac{F(dy)}{y^2}
\label{intro:0}
\end{equation}
where $\Lambda = F$, a probability measure. A Wright-Fisher generator equation, identical to (\ref{intro:0}) is
\begin{equation}
{\cal L}g(x) = \frac{1}{2}x(1-x)\mathbb{E}\Big [g^{\prime\prime}\big (x(1-W)+VW\big )\Big ]
\label{intro:2}
\end{equation}
where $W=UY$, $Y$ has distribution $F$, $U$ has a density $2u$, $u \in (0,1)$, $V$ is uniform on $(0,1)$, and $U,V,Y$ are independent. If $W=0$ the usual Wright-Fisher generator is obtained. The equation (\ref{intro:2}) is very suggestive of a strong representation between the $\Lambda$-Fleming-Viot and Wright-Fisher processes.

The $d$-dimensional generator has polynomial eigenvectors and eigenvalues which are analogues of those in the Wright-Fisher generator. The eigenvalues are
\[
\frac{1}{2}n(n-1)\mathbb{E}\Big[(1-W)^{n-2}\Big ],\>n=2,3,\ldots
\]
which are equal to the $\Lambda$-coalescent total merger rates from $n$ blocks. If $d=2$ the polynomial eigenvectors are analogues of the Jacobi polynomials. 

The two-dimensional process is considered in detail in this paper. An integral equation is found for the stationary distribution when  there is mutation. This leads to an interesting equation for the frequency spectrum in the infinitely-many-alleles $\Lambda$-Fleming-Viot model when the $\Lambda$-coalescent comes down from infinity. If frequencies of the alleles are denoted by $x_{(1)} \geq x_{(2)} \geq \cdots $ and $\mathbb{E}$ denotes expectation in the stationary distribution then the (1-dimensional) frequency spectrum $\beta(x)$ is defined by
\begin{equation}
 \mathbb{E}\Big [\sum_{k=1}^\infty f(x_{(k)})\Big ]
 = \int_0^1f(x)\beta(x)dx
 \label{fsdef:0}
\end{equation}
where $f \in C([0,1])$ and $f(x)/x$ is bounded as $x \to 0$. 
The 1-dimensional frequency spectrum is the same as the first factorial moment measure for the allele frequencies $\{x_{(i)}\}$ regarded as a point process. Equation (\ref{fsdef:0}) follows from general point processes theory \citep{DVJ2005}.
From the definition (\ref{fsdef:0}) it follows that $z\beta(z)$, $0 < z < 1$ is a probability density.  Let $Z$ be a random variable with this density, $Z_*$ a random variable size-biassed with respect to $Z$, $Z^*$ a random variable size-biassed with respect to $1-Z$ and $V$ a uniform random variable on $[0,1]$. Then
\begin{equation}
VZ_* =^{\cal D} (1-W)Z^* + VW,
\label{intro:3}
\end{equation}
where the random variables are independent of each other. The left side is the limit distribution of excess life in a renewal process with increments distributed as $Z$ \citep{C1970}, so the equation suggests a renewal process.  We do not have a probabilistic solution of (\ref{intro:3}) which would possibly lead to knowing $\beta(z)$.

In a two-dimensional process with no mutation and genic selection a proof is given that there is a critical selection value such that if the selection coefficient is greater or equal to the critical value then fixation, when the boundary 1 is hit, has probability 1 beginning from any non-zero frequency. This is an analytic proof different from proofs by \cite{DEP2011} and \citet{F2013} which uses our particular representation of the generator. A computational solution for the probability of fixation, when fixation is not certain, is found which is analogous to that in the Wright-Fisher model. \citet{BP2013} construct a lookdown process \citep{DK1996} in this model.

The moment dual process in the Fleming-Viot process is the usual $\Lambda$-coalescent back in time. In a model with two types, generator (\ref{intro:0}), and $X(t)$ the frequency of the first type at time $t$ there is a dual equation
\[
\mathbb{E}_{X(0)=x}\big [X(t)^n\big ]
= \mathbb{E}_{L(0)=n}\big [x^{L(t)}\big ].
\]
In this equation $\{L(t)\}_{t \geq 0}$ is a $\Lambda$-coalescent process back in time with transition rates $\lambda_{nk}$. Expectation on the left is with respect to $X(t)$, and on the right with respect to $L(t)$.

 In the Wright-Fisher representation using a different set of polynomials $g_n(x)$ which mimic $x^n$ in the usual Wright-Fisher diffusion as test functions produces a dual death process which has a similarity to the Kingman coalescent and decreases by units of one. The $d$-dimensional version $g_{\bm{n}}({\bm x})$ analogous to $\bm{x}^{\bm n}$ has an expectation in the stationary distribution of a model with parent independent mutation that is similar to a Dirichlet moment
\begin{equation*}
\mathbb{E}\big [g_{\bm{n}}(\bm{x}) \big ]
= \frac{
\prod_{i=1}^d\Bigl [\prod_{j=1}^{n_i}\bigl ((j-1)\mathbb{E}\bigl [(1-W)^{j-2}\bigr ]+\theta_i\bigr )\Bigr ]
}
{
\prod_{j=1}^{n}\bigl ((j-1)\mathbb{E}\bigl [(1-W)^{j-2}\bigr ]+\theta\bigr )
}.
\end{equation*}
Bold face notation will be used for $d$-dimension vectors in the paper, and the shorthand notation
 $\bm{x}^{\bm{n}}\equiv \prod_{i=1}^dx_i^{n_i}$.
There is even an analogue of the Ewens' sampling formula in the Poisson Dirichlet process of
\[
\frac{n!\theta^k}{n_1\cdots n_k}
\cdot \frac{\prod_{i=1}^{k}\Bigl [\prod_{j=2}^{n_i}\mathbb{E}\big [(1-W)^{j-2}\big ]\Bigr ]}
{
\prod_{j=1}^{n}\Big [ (j-1)\mathbb{E}\big [(1-W)^{j-2}\big ]+\theta\Big ]
}.
\]

There are many intriguing analogues between the $\Lambda$-Fleming-Viot process and the Wright-Fisher diffusion process which come from the generator representation. 

Exact calculations are always likely to be difficult because of the jump process nature of the $\Lambda$-Fleming-Viot process. A first step in this direction, for certain classes of Fleming-Viot processes where stationary distributions are characterized, can be found in \citet{H2012}.

%%%%
\section{A Wright-Fisher generator connection}
The $\Lambda$-Fleming-Viot  generator has an interesting connection with a Wright-Fisher diffusion generator that we now develop.
\medskip

%% Theorem
\noindent
{\bf Theorem 1.}
\emph{ 
Let ${\cal L}$ be the  $\Lambda$-Fleming-Viot  generator (\ref{L:basic}), $V$ be a uniform random variable on $[0,1]$, $U$ a random variable on $[0,1]$ with density $2u,\>0< u <1$ and $W=YU$, where 
$Y$ has distribution $F$ and $V,U,Y$ are independent. Denote the first and second derivatives of a function $g(\bm{x})$ in ${\cal C}^2(\Delta)$ by
\[
g_i(\bm{x}) = \frac{\partial}{\partial x_i}g(\bm{x}),\>g_{ij}(\bm{x}) = \frac{\partial^2}{\partial x_i\partial x_j}g(\bm{x}).
\]
Then
\begin{equation}
{\cal L}g(\bm{x}) = \frac{1}{2}\sum_{i,j=1}^dx_i(\delta_{ij}-x_j)
\mathbb{E}\Bigl [g_{ij}\big (\bm{x}(1-W)+WV\bm{e}_i\big )\Bigr ],
\label{L:LWF}
\end{equation}
where expectation $\mathbb{E}$ is taken over $V,W$. 
}
\medskip

\noindent
{\bf Proof.}
Taking the expectation with respect to $V$ the right side of  (\ref{L:LWF}) is equal to 
\begin{eqnarray}
&&\frac{1}{2}\sum_{i,j=1}^dx_i(\delta_{ij}-x_j)
\mathbb{E}\Biggl [\frac{g_j\big (\bm{x}(1-W)+W\bm{e}_i\big )
-g_{j}\big (\bm{x}(1-W)\big )}{W}\Biggr ]
\nonumber \\
&&~= 
\int_0^1\Bigg [\sum_{i,j=1}^dx_i(\delta_{ij}-x_j)
\int_0^1g_j\big (\bm{x}(1-uy)+uy\bm{e}_i\big )du\Bigg ]\frac{F(dy)}{y}
\nonumber \\
&&~~~~~
-\int_0^1\Bigg [\sum_{i,j=1}^dx_i(\delta_{ij}-x_j)\int_0^1g_j\big (\bm{x}(1-uy)\big )du\Bigg ]\frac{F(dy)}{y}.
\label{L:LWF:a}
\end{eqnarray}
To simplify (\ref{L:LWF:a}) note that 
\begin{eqnarray*}
\frac{\partial}{\partial u}g\big (\bm{x}(1-uy)\big )
&=& -y\sum_{j=1}^dx_jg_j\big (\bm{x}(1-uy)\big )\>\>\text{and}
\nonumber \\
\frac{\partial}{\partial u}g\big (\bm{x}(1-uy)+uy\bm{e}_i\big )
&=& y\sum_{j=1}^d(\delta_{ij}-x_j)g_j\big (\bm{x}(1-uy)+uy\bm{e}_i\big ).
\end{eqnarray*}
Therefore (\ref{L:LWF:a}) is equal to 
\begin{eqnarray}
&&~\int_0^1\Bigg [\sum_{i=1}^dx_i\int_0^1\frac{\partial}{\partial u}g\big (\bm{x}(1-uy)+uy\bm{e}_i\big )du\Bigg ]\frac{F(dy)}{y^2}
\label{intermediate:1}
\\
&&- \int_0^1\Bigg [\int_0^1\Big (1-\sum_{i=1}^dx_i\Big )
\frac{\partial}{\partial u}g\big (\bm{x}(1-uy)\big )du\Bigg ]\frac{F(dy)}{y^2}
\label{intermediate:2}
\\
&&~~= \int_0^1\Bigg [\sum_{i=1}^dx_i\big (g(\bm{x}(1-y)+y\bm{e}_i) - g(\bm{x})\big )\Bigg ]\frac{F(dy)}{y^2}.
\label{intermediate:3}
\end{eqnarray}
In the calculation the term (\ref{intermediate:1}) is equal to (\ref{intermediate:3}), and the term (\ref{intermediate:2}) vanishes.
\qed
\medskip

The Wright-Fisher generator (\ref{L:WF}) is included in (\ref{L:LWF}) when $W\equiv 0$.
\medskip

\noindent
{\bf Corollary 1.}
\emph{
$\{X_1(t)\}_{t\geq 0}$ is a Markov process with generator acting on functions in ${\cal C}^2([0,1])$ specified by
\begin{equation}
{\cal L}g(x) = \frac{1}{2}x(1-x)\mathbb{E}\Bigl [g^{\prime\prime}\big (x(1-W) + WV\big )\Bigr ].
\label{L:WFM}
\end{equation}
}
\medskip

\noindent
{\bf Proof.} Let $g(\bm{x})$ in (\ref{L:LWF}) be a function of the first co-ordinate only, then (\ref{L:WFM}) follows easily, with $x\equiv x_1$.
\qed
\medskip

$W$ possibly has an atom at $0$,
$P(W=0) = P(Y=0)$,
and is continuous for $W >0$ with a density
\begin{equation}
f_W(w)=2wF^+(w),
\label{WY:1}
\end{equation}
where
\begin{equation}
F^+(w) = \int_w^1y^{-2}F(dy).
\label{F+def}
\end{equation}
There is a correspondence between $F$ and the distribution of $W$. Given a random variable $W$ with a possible atom at $0$ and a density $f_W(w)$, $0 < w \leq 1$, then there exists independent random variables $U,Y$, where $U$ has density $2u$, $0 < u < 1$ such that $W=YU$ if and only if $f_W(1)=0$ and $f_W(w)/w$ is decreasing in $(0,1]$. Possible densities for the continuous component of $W$ are proportional to the Beta $(a,b)$ densities with $a \leq 2$ and $b \geq 1$. In particular if $Y$ has a Beta $(2-\alpha,\alpha)$ distribution, then $W$ has a Beta $(2-\alpha,1+\alpha)$ distribution.
%Reallocate thm
\medskip

The next theorem gives a connection between $W$, the $\Lambda$-coalescent rates and the Laplace exponent.
\medskip

\noindent
{\bf Theorem 2.}
\emph{
%Let $n\geq 2$
\begin{eqnarray}
\sum_{k=2}^n{n\choose k}\lambda_{nk}
&=&
\int_0^1
\Bigl [1 - (1-y)^n - ny(1-y)^{n-1}\Bigr ]\frac{F(dy)}{y^2}
\nonumber \\
&=& \frac{1}{2}n(n-1)\mathbb{E}\Bigl [(1-W)^{n-2}\Bigr ],\>\text{for~~}n \geq 2. 
\label{connection:e}
\end{eqnarray}
The individual rates (\ref{lnk:rates}) can be expressed for
$2 \leq k \leq n$ as
\begin{eqnarray}
{n\choose k}\lambda_{nk} &= &{n\choose k}\int_0^1y^k(1-y)^{n-k}\frac{F(dy)}{y^2}
\nonumber \\
&=& \frac{n}{2}\mathbb{E}\Big [P_{k-1}(n,W) - P_{k}(n,W)\Big ],
\label{negbin:0}
\end{eqnarray}
where
\[
P_k(n,w) = {n-1\choose k}(1-w)^{n-k-1}w^{k-1}.
\]
The Laplace exponent
\begin{equation}
\psi(q) = \frac{q}{2}\mathbb{E}\Bigg [\frac{1-e^{-qW}}{W}\Bigg ].
\label{LFconnection:0}
\end{equation}
}
\medskip

\noindent
{\bf Proof.}
\begin{eqnarray*}
\frac{1}{2}n(n-1)\mathbb{E}\Bigl [(1-W)^{n-2}\Bigr ]
&=& 
\frac{1}{2}n(n-1)\int_0^1\int_0^1(1-uy)^{n-2}2uduF(dy)
\nonumber \\
&=&
\int_0^1\int_0^1u\frac{\partial^2}{\partial u^2}(1-uy)^{n}du\frac{F(dy)}{y^2}
\nonumber \\
&=&
\int_0^1
\Bigl [1 - (1-y)^n - ny(1-y)^{n-1}\Bigr ]\frac{F(dy)}{y^2}
\nonumber \\
&=&
\sum_{k=2}^n{n\choose k}\lambda_{nk}.
\end{eqnarray*}
For the individual rates, showing (\ref{negbin:0}) is an exercise in integration by parts which follows from
\begin{eqnarray*}
&&\frac{n}{2}\mathbb{E}\big [P_k(n,W)\big ]
\\
&&~=\frac{n}{2}{n-1\choose k}\int_0^1\int_0^1(1-uy)^{n-k-1}(uy)^{k-1}2uduF(dy)
\\
&&~=-n{n-1\choose k}(n-k)^{-1}\int_0^1y^{k-2}\int_0^1u^k\frac{\partial (1-uy)^{n-k}}{\partial u}duF(dy)
\\
&&~= -\textcolor{black}{{n\choose k}}\lambda_{nk} + \frac{n}{2}\mathbb{E}\big [P_{k-1}(n,W)\big ].
\end{eqnarray*}
Note that $P_n(n,w) \equiv 0$ so
$
\lambda_{nn} = \frac{n}{2}\mathbb{E}\big [P_{n-1}(n,W)\big ]
$.
To show (\ref{LFconnection:0}) 
\begin{eqnarray*}
\psi(q) &=& 
\int_0^1 \big (e^{-qy} -1 +qy\big )\frac{F(dy)}{y^2}
\nonumber \\
&=&\int_0^1\sum_{k=2}^\infty (-1)^k\frac{q^ky^{k-2}}{k!}F(dy)
\\
&=& 
\frac{1}{2}\int_0^1\sum_{k=2}^\infty (-1)^k\frac{q^k}{(k-1)!}
\int_0^1(uy)^{k-2}2uduF(dy)
\\
&=&
 \frac{q}{2}\mathbb{E}\Bigg [\frac{1-e^{-qW}}{{\small W}}\Bigg ]
\end{eqnarray*}
\qed
\medskip

%The Laplace transform of $W$is related to the Laplace exponent (\ref{Laplaceexponent}) by
%\begin{equation}
%\mathbb{E}\bigl [e^{-qW}\bigr ] = 2q^{-2}\psi (q).
%\label{connection:0}
%\end{equation}
%The condition (\ref{condition:0}) of \citet{BLG2006} for the $\Lambda$-coalescent to come down from infinity is equivalent to 
%\begin{equation}
%\int_1^\infty
%\frac{dq}{q^2\mathbb{E}\big [e^{-qW}\big ]} < \infty
%\label{condition:3}
%\end{equation}
%because of the connection (\ref{condition:0}).

The random variables $Y,W,V$ from Theorem 1 are used frequently in the paper, so their definition will be assumed.

%%%%%
\subsection{Mutation and selection} Mutation can be added to the model by assuming that mutations occur at rate $\theta/2$ and changes of type $i$ to type $j$ are made according to a transition matrix $P$. This is equivalent to mutations occurring at rate $\theta/2$ on the dual $\Lambda$-coalescent tree. The generator (\ref{L:basic}) then has an additional term added of
\begin{equation}
\frac{\theta}{2}\sum_{i=1}^d\Big (\sum_{j=1}^dp_{ji}x_j - x_i\Big )\frac{\partial}{\partial x_i}.
\label{mutation:pij}
\end{equation}
If mutation is parent independent $\theta p_{ji}= \theta_i$, not depending on $j$, and the additional term simplifies to
\begin{equation}
\frac{1}{2}\sum_{i=1}^d\Big (\theta_i - \theta x_i\Big )\frac{\partial}{\partial x_i}.
\label{mutation:pim}
\end{equation}
If $d=2$ and $x_1=x$, $x_2=1-x$, then the generator acting on functions $g(x)$ in ${\cal C}^2([0,1])$ is specified by
\begin{eqnarray}
{\cal L}g\big (x\big ) &=& 
\int_0^1 \Bigg [x\Big (g\big (x(1-y)+y\big ) - g\big (x\big )\Big )
\nonumber \\
&&~~+ (1-x)\Big (g\big (x(1-y)\big ) - g\big (x\big )\Big )
\Bigg ]
\frac{F(dy)}{y^2} 
\nonumber \\
&&~~
+ \textcolor{black}{\frac{1}{2}}(\theta_1 - \theta x)g^{\prime}(x).
\label{two:00}
\end{eqnarray}
\citet{H2012} finds the stationary distribution in a process with generator specified by
\begin{eqnarray}
{\cal L}_Hg\big (x\big ) &=&
\int_0^1 \Big [x\Big (g\big (x(1-y)+y\big ) - g\big (x\big )\Big )
\nonumber \\
&&+ (1-x)\Big (g\big (x(1-y)\big ) - g\big (x\big )\Big )\Big ]\frac{B_{1-\alpha,1+\alpha}(dy)}{y^2} 
\nonumber \\
&&+\int_0^1 \Big [\theta_1 g\big (x(1-y)+y\big ) 
+ \theta_2g\big (x(1-y)\big ) - \theta g\big (x\big )\Big ]\frac{B_{1-\alpha,\alpha}(dy)}{(\alpha+1)y},
\nonumber \\
\label{two:H}
\end{eqnarray}
where $0 < \alpha < 1$ and $B_{a,b}(dy)$ denotes a Beta $(a,b)$ density. In his model there is simultaneous mutation, where at rate
$\theta_1B_{1-\alpha,\alpha}(dy)/\Big ((\alpha+1)y \Big )$ a proportion $y$ of the population is replaced by type 1 individuals and similarly at rate $\theta_2B_{1-\alpha,\alpha}(dy)/\Big ((\alpha+1)y \Big )$ a proportion $y$ of the population is replaced by type 2 individuals. This is an unusual mutation mechanism and the generators  (\ref{two:00}) and  (\ref{two:H}) are different even when $F=B_{1-\alpha,1+\alpha}$.

\citet{EGT2010} study a $\Lambda$-Fleming-Viot process with viability selection whose generator acting on functions in ${\cal C}^2(\Delta)$ takes the form
\begin{eqnarray}
{\cal L}g(x)
&=& 
\int_0^1 \sum_{i=1}^d x_i\big (g(\bm{x}(1-y)+ y\bm{e}_i) - g(\bm{x})\big )\frac{F(dy)}{y^2}
\nonumber \\
&&~~~~~~~~
-\int_0^1 \sum_{i=1}^d x_i\big (g(\bm{x}(1-y)+ y\bm{e}_i) - g(\bm{x})\big )\frac{K_i(dy) }{y}
\nonumber \\
&&~~~~~~~~
+\frac{\theta}{2}\sum_{i=1}^d\Big (\sum_{j=1}^dp_{ji}x_j - x_i\Big )\frac{\partial}{\partial x_i}
g(\bm{x}).
\label{Lambdagen:1}
\end{eqnarray}
To describe the measures in (\ref{Lambdagen:1}) let $G_i$, $i\in [d]$ be the $\Lambda$-measures for the individual types, which are positive measures on $[0,1]$ and $F$ be a reference measure such that 
\[
K_i(dy) = \frac{F(dy) - G_i(dy)}{y} 
\]
are bounded signed measures on $[0,1]$. A selection model analogous to the Wright-Fisher model with genic selection (see for example \citet{E2004}) is obtained by taking  
\[
K_i(\cdot) = \sigma_i\delta_\epsilon (\cdot),
\]
and letting $\epsilon \to 0+$. Selection is very weak 
in this limit in the sense that a limit is taken where all the measures approach $F$, whereas there is a much larger effect when the measures $G_i$ are different.
The corresponding sequence of generators converges to
\begin{eqnarray}
{\cal L}^{\bm{\sigma}} g(\bm{x}) 
&=& 
\int_0^1 \sum_{i=1}^d x_i\big (g(\bm{x}(1-y)+ y\bm{e}_i) - g(\bm{x})\big )\frac{F(dy)}{y^2}
\nonumber \\
&&~~~~~~~~
-\sum_{i=1}^dx_i \big (\sigma_i - \sum_{k=1}^d\sigma_k x_k\big )
\frac{\partial}{\partial x_i}g(\bm{x})
\nonumber \\
&&~~~~~~~~
+\frac{\theta}{2}\sum_{i=1}^d\big (\sum_{j=1}^dp_{ji}x_j - x_i\big )\frac{\partial}{\partial x_i}
g(\bm{x}).
\label{Lambdagen:3}
\end{eqnarray}
\citet{EGT2010} find the dual Lambda coalescent corresponding to (\ref{Lambdagen:1}) and (\ref{Lambdagen:3}).
%%%
%%Replaced section
%%

\subsection*{\textcolor{black}{Fixation probability with selection when $d=2$ types}}
If there are $d=2$ types, no mutation, $X=X_1$, $\sigma_1\leq 0$, $\sigma_2=0$, then with notation $\beta = -\sigma_1 \geq 0$ the generator equation (\ref{Lambdagen:3}) reduces to
\begin{equation*}
{\cal L}^\beta g(x) = \frac{1}{2}x(1-x)\mathbb{E}\Bigl [g^{\prime\prime}(x(1-W)+WV)\Bigr ]
+\beta x(1-x)g^\prime(x)
%\label{twosel:0}
\end{equation*}
Let $P(x)$ be the probability that the first type fixes, starting from an initial frequency of $x$. Then $P(0)=0$, $P(1)=1$, and $P(x)$ is the solution of
\[
{\cal L}^\beta P(x) = 0.
\]
That is
\begin{equation}
\mathbb{E}\Bigl [P^{\prime\prime}(x(1-W)+WV)\Bigr ]
+2\beta P^\prime(x)=0,
\label{fix}
\end{equation}
and taking the expectation with respect to $V$, 
\begin{equation}
\mathbb{E}\Biggl [\frac{P^{\prime}\big(x(1-W)+W\big) - P^\prime\big (x(1-W)\big)}{W}\Biggr ]
+2\beta P^\prime(x)=0.
\label{fix:1}
\end{equation}
Integrating and taking care of a possible discontinuity $P(0+)$  at $x=0$,
\begin{eqnarray}
&&\mathbb{E}\Biggl [\frac{P\big(x(1-W)+W\big) - P\big (x(1-W)\big) - P(W) + P(0+)}{W(1-W)}\Biggr ]
\nonumber \\
&&~~~~~~~~+2\beta \big [P(x)-P(0+)\big ]
=0.
\label{fix:1a}
\end{eqnarray}
Alison Etheridge and Jay Taylor have obtained equivalent formulae to (\ref{fix:1},\ref{fix:1a}) in the Beta coalescent using integration by parts, private communication (2008).
\citet{DEP2011,DEP2012} study fixation probabilities in the $\Lambda$-coalescent. An interesting feature is that for some $\Lambda$-measures and $\beta$ it can happen that $P(x) = 1$ or $P(x) = 0$ for all $x \in (0,1)$. They show that fixation is certain (that is, $P(x)=1$, $x \in (0,1]$) if and only if
\begin{equation}
\beta \geq \beta^* = - \int_0^1\frac{\log (1-y)}{y^2}F(dy)
\label{condition:0}
\end{equation}
under the assumption that $\beta^* < \infty$.
% Bob: Deleted in version 5: equivalent to $\int_0^1y^{-1}F(dy) < \infty$. 
%
If $\beta^* = \infty$ then fixation is not certain.
Their proof is for the Eldon-Wakeley coalescent where $F$ has a single point of increase in $(0,1]$.
The general formula (\ref{condition:0}) is mentioned in the paper and has an analogous proof to the Eldon-Wakeley case,  private communication (2013).
They use a clever comparision of $P(x)$ with sub-harmonic and super-harmonic functions. If $u(x)$ is such that $u(0)=0$, $u(1)=1$ then if
${\cal L}^\beta u(x) \leq 0$ for all $x \in (0,1)$ they show that $P(x) \leq u(x)$ for all $x \in (0,1)$. Similarly if ${\cal L}^\beta u(x) \geq 0$ for all $x \in (0,1)$, $P(x) \geq u(x)$ for all $x \in (0,1)$. Comparison functions used are $u(x) = x^p$, and $u(x) = Cx^p + (1-C)x$, $ 0 < p < 1$ and $C > 1$. \citet{F2013} gives an elegant martingale proof based on a dual process that (\ref{condition:0}) is necessary and sufficient for $P(x)=1$, $x \in (0,1]$, though does not include the critical case when $\beta = \beta^*$ in his proof.
Another way to express (\ref{condition:0}) is
\begin{equation}
2\beta \geq 2\beta^* = \mathbb{E}\Big [\frac{1}{W(1-W)}\Big ].
\label{condition:1}
\end{equation}
%
% Bob: Deleted in version 5: with $\int_0^1y^{-1}F(dy) < \infty$ equivalent to $\mathbb{E}\big [W^{-1}\big ] < \infty$. The equivalence is straightforward to show:
\begin{eqnarray*}
\frac{1}{2}\mathbb{E}\Big [\frac{1}{W(1-W)}\Big ] &=&
\frac{1}{2}\int_0^1\int_0^1\frac{1}{uy(1-uy)}2udu F(dy)
\nonumber \\
&=& 
\int_0^1\int_0^1\frac{du}{1-uy}\frac{F(dy)}{y}
\nonumber \\
&=& 
\int_0^1\frac{-\log (1-y)}{y^2}F(dy).
\end{eqnarray*}
For interest we show how our representation can be used to give a proof when $\beta^* < \infty$.
\medskip

\noindent
{\bf Theorem 3.}\emph{ \citep{DEP2011,DEP2012,F2013}.
Let $\beta^* < \infty$. Then $P(x) = 1$ for all $x \in (0,1]$ if and only $\beta \geq \beta^*$.%  If $\beta^* = \infty$ then  $P(x) < 1$ for all $x \in [0,1)$.
}
\medskip

\noindent
{\bf Proof.}
\medskip

\noindent\emph{If.} Let $\beta=\beta^*$. For $x \in (0,1]$, from (\ref{fix:1a})
\begin{equation}
0=\mathbb{E}\Biggl [\frac{P\big(x(1-W)+W\big) - P\big (x(1-W)\big) - P(W) + P(x)}{W(1-W)}\Biggr ]
\label{delicate:0}
\end{equation}
$P(x)$ is a non-decreasing function of $x$ and
since the right side of (\ref{delicate:0}) must be zero, with probability 1,
\[
P(x(1-W)+W) - P(W) = 0\text{~~and~~}
P(x) - P(x(1-W))  = 0.
\]
This can only be true if $P(x)=1$ for all $x \in (0,1]$, since $P(1)=1$.
Now take $\beta \geq \beta^*$. $P_\beta(x)\equiv P(x)$ is a non-decreasing function of $\beta$ for fixed $x$ because a higher selective parameter produces a higher probability of fixation. Thus $P_\beta(x) \geq P_{\beta^*}(x) = 1$ for all $x \in (0,1]$ and it must be that $P_\beta(x) = 1$.
\medskip

\noindent\emph{Only if.} Let $\beta < \beta^*< \infty$ and suppose that $P(x)=1$ for $x\in (0,1]$. We show this assumption is contradictory. Consider a test function 
\[
v(x) = \log(x) + K(1-x),
\]
where $K>0$ is a constant.
A generator equation is that
\begin{equation}
\mathbb{E}_x\big [v(X(t))\big ] - v(x) = \int_0^t\mathbb{E}_x\big [{\cal L}^\beta v(X(u))\big ]du.
\label{newth:0}
\end{equation}
Equation (\ref{newth:0})  evaluates to 
\begin{equation}
\mathbb{E}_x\big [v(X(t))\big ] - v(x)
= \frac{1}{2}\int_0^t \mathbb{E}_x\Big [ \big (1-X(u) \big )A(u)\Big ]du,
\label{newth:1}
\end{equation}
where
\begin{eqnarray}
A(u) &=&
\mathbb{E}\Bigg [ X(u)\frac {
\big (X(u)(1-W) + W\big )^{-1} - \big (X(u)(1-W)\big )^{-1}
}
{W}
\nonumber \\
&&~~~~~~~~~~+ X(u)2\beta X(u)^{-1} - 2K\beta X(u)\Bigg ]
\nonumber \\
&=& \mathbb{E}\Bigg [
-\frac{1}{\Big (X(u)(1-W)+W\Big )(1-W)} - 2K\beta X(u)\Bigg ] + 2\beta.
\label{newth:2}
\end{eqnarray}
Choose $K$ large enough so that the minimum value over $x \in [0,1]$ of
\begin{equation}
\mathbb{E}\Bigg [
\frac{1}{\Big (x(1-W)+W\Big )(1-W)} + 2K \beta x\Bigg ]
\label{choose:K}
\end{equation}
is attained when $x=0$. Then $A(u) \leq -2\beta^* + 2\beta < 0$.
Let $t \to \infty$ in (\ref{newth:1}). $X(t) \to 1$ with probability 1, so $\mathbb{E}_x\big [\log\big (X(t)\big ) + K\big (1-X(t)\big ] \to 0$
 and the limit equation is
\begin{eqnarray}
-\log x - K(1-x) &=&  \frac{1}{2}\int_0^\infty \mathbb{E}_x\Big [ (1-X(u))A(u)\Big ]du
\nonumber \\
&\leq& (\beta - \beta^*)\int_0^\infty \mathbb{E}_x\Big [ (1-X(u))\Big ]du < 0.
\label{newth:6}
\end{eqnarray}
Choose $x$ small enough so that the left side of (\ref{newth:6}) is positive. Then the signs of both sides of (\ref{newth:6}) are contradictory.
Therefore the assumption that $P(x)=1$ for all $x\in (0,1]$ is contradictory. Let $x_0$ be the maximal point where $P(x) < 1$ for 
$0 < x < 1$. It cannot happen that $P(x) = 1$ for $x_0 \leq x < 1$. Suppose that this does occur. Let $X(0)=x_0$ and consider the local exit behaviour of $X$ in $(x_0-\epsilon,x_0+\epsilon)$, for small $\epsilon > 0$. For small enough $\epsilon$ there is positive probability that there is a path where $X$ first exits the interval at less than or equal to $x_0 - \epsilon$. The strong Markov propery of $X$ then implies $P(x_0-\epsilon) = 1$, which is contradictory.
Therefore $P(x) < 1$ for all $x \in [0,1)$.
\qed
%%%
\medskip

In the Kingman coalescent (\ref{fix}) becomes
\[
P^{\prime\prime}(x)  + 2\beta P^\prime(x) = 0,
\]
with a solution
\begin{equation}
P(x) = \frac{1 - e^{-2\beta x}}{1-e^{-2\beta}}.
\label{Kfix:0}
\end{equation}
We provide a computational solution for $P(x)$ in the $\Lambda$-Fleming-Viot model when fixation or loss is not certain from $x \in (0,1)$ that imitates (\ref{Kfix:0}). A sequence of polynomials $\{h_n(x)\}_{n=0}^\infty$ that is used in the proof is defined as the solutions of
\begin{equation}
\mathbb{E}\Biggl [\frac{h_n\big(x(1-W)+W\big) - h_n\big (x(1-W)\big)}{W}\Biggr ]
= nh_{n-1}(x),
\label{polysel:1}
\end{equation}
where the leading coefficient in $h_n(x)$ is 
\begin{equation}
\frac{1}{\prod_{j=1}^{n-1}\mathbb{E}\bigl [(1-W)^j\bigr ]}.
\label{polysel:1a}
\end{equation}
This choice makes the coefficients of $x^{n-1}$ in (\ref{polysel:1}) agree.
The argument in the expectation in (\ref{polysel:1}) is interpreted as $h_n^\prime(x)$ at $W=0$. There is a family of polynomial solutions to (\ref{polysel:1}) depending on an arbitary recursive choice of constant coefficients.
The constant coefficients in the polynomials are chosen carefully to obtain a solution for the fixation probability. The polynomials $h_n(x)$ imitate $x^n$ and are equal if $W\equiv 0$. Let $h_0(x) = 1$, and 
\[
h_n(x)=\sum_{r=0}^na_{nr}x^r.
\]
Then from (\ref{polysel:1}) for $j=n-2,\ldots ,0$
\[
\sum_{j=0}^{n-1}\sum_{r=j+1}^n{r \choose j}\mathbb{E}\Bigl [(1-W)^jW^{r-j-1}\Bigr ]a_{nr}x^j
= n\sum_{j=0}^{n-1}a_{n-1j}x^j,
\]
so equating coefficients of $x^j$ on both sides,
\begin{equation}
\sum_{r=j+1}^n{r \choose j}\mathbb{E}\Bigl [(1-W)^jW^{r-j-1}\Bigr ]a_{nr} = na_{n-1j}.
\label{polysel:1b}
\end{equation}
Given the coefficients $\{a_{n-1j}\}_{j=0}^{n-1}$ of $h_{n-1}(x)$ the coefficients of $h_n(x)$, $\{a_{nj}\}_{j=1}^n$ are recursively determined by choosing $a_{nn}$ from (\ref{polysel:1a}), then taking $j=n-1,\ldots ,0$ in (\ref{polysel:1b}).
There is an arbitrary choice of $a_{n0}$ that needs to be made at this stage to progress with the recursion. 
\medskip

\noindent
{\bf Theorem 4.}
\emph{
Let $0< \beta < \beta^*$. The fixation probability
\begin{equation*}
P(x) = \bigl (1 - e^{-2\beta} \bigr )^{-1}\sum_{n=1}^\infty (-1)^{n-1}\frac{(2\beta)^n}{n!}H_n(x),
%\label{polysel:2}
\end{equation*}
where $\{H_n(x)\}$ are polynomials derived from
\[
H_n(x) = \int_0^x nh_{n-1}(\xi)d\xi
\]
with the constants $\{h_n(0)\}$ chosen so that 
\begin{equation}
\int_0^1nh_{n-1}(\xi)d\xi = 1.
\label{choice}
\end{equation}
}
\medskip

\noindent
{\bf Proof.}
\medskip
Try a series solution
\begin{equation}
P^\prime(x) =
B(\beta)\sum_{n=1}^\infty (-1)^{n-1} (2\beta)^nc_nh_{n-1}(x),
\label{polysel:3}
\end{equation}
where $\{h_n(x)\}$ satisfies (\ref{polysel:1}), $B(\beta)$ is a constant, and $\{c_n\}$ are constants not depending on $\beta$.  Then substituting in (\ref{fix:1}) 
\[
\sum_{n=2}^\infty (-1)^{n-1}(2\beta)^nc_n(n-1)h_{n-2}(x) + (2\beta)\sum_{n=1}^\infty(-1)^{n-1}(2\beta)^nc_nh_{n-1}(x) = 0.
\]
This identity is satisfied if $c_1=-1$, without loss of generality, and
\[
c_n=-\frac{1}{(n-1)!},\>n=2,3,\ldots
\]
Integrating in (\ref{polysel:3})
\[
P(x) =
B(\beta)\sum_{n=1}^\infty (-1)^{n-1}\frac{(2\beta)^n}{n!}\int_0^xnh_{n-1}(\xi)d\xi.
\]
Choosing (\ref{choice}) to hold and knowing $P(1)=1$ shows that 
\[
B(\beta) = \bigl (1 - e^{-2\beta} \bigr )^{-1}.
\]
\qed
\medskip

\noindent
{\bf Corollary 2.} \emph{
A computational solution for $P(x)$ is found from evaluating the polynomials
\[
H_n(x) = \sum_{r=1}^nb_{nr}x^r,
\]
where $H_1(x)=x$ and the coefficients $\{b_{nr}\}$ are defined recursively from
\begin{equation}
\sum_{r=j}^n{r\choose j-1}\mathbb{E}\Bigl [(1-W)^{j-1}W^{r-j}\Bigr ](r+1)b_{n+1r+1}
= (n+1)jb_{nj},
\label{poly:recurse}
\end{equation} 
with 
\[
b_{n+11} = 1 - \sum_{j=2}^{n+1}b_{n+1j}
\]
for $n=1,2, \ldots$, $j=n-1,\ldots ,1$.
}
(\ref{poly:recurse}) is equivalent to
\begin{equation}
2\sum_{r=j+1}^{\textcolor{black}{n+1}}\Bigg [\sum_{k=r-j+1}^{r}{r \choose k}\lambda_{rk}\Bigg ]b_{n+1r}
= (n+1)jb_{nj}.
\label{poly:recurse:0}
\end{equation} 
\medskip

\noindent
{\bf Proof.} Relating the coefficients of $H_n(x)$ to those of $h_{n-1}(x)$ 
\[
b_{nj} = \frac{n}{j}a_{n-1j-1},\>j=2,\ldots ,n\text{~~~and~~}b_{n1} = 1 - \sum_{j=2}^nb_{nj}.
\]
Substituting in (\ref{polysel:1b}) and shifting the index $j \to j+1$ completes the proof of (\ref{poly:recurse}). The alternative form (\ref{poly:recurse:0}) is found by noting that
\begin{eqnarray*}
\frac{r+1}{2}\mathbb{E}\big [P_{\textcolor{black}{r-j+1}}(r+1,W)\big ]
&= &\frac{r+1}{2}{r \choose j-1}
\mathbb{E}\big [(1-W)^{j-1}W^{r-j}\big ]
\\
&=&
\sum_{\textcolor{black}{k=r-j+2}}^{r+1}{r+1 \choose k}\lambda_{r+1k}
\end{eqnarray*}
from (\ref{negbin:0}),
substituting, then shifting the index of summation $r\to r+1$.
\qed
\subsection{Eigenstructure of the $\Lambda$-Fleming-Viot process}
The generator of the $\Lambda$-Fleming-Viot process (\ref{L:LWF}) with mutation term (\ref{mutation:pij})
\begin{eqnarray}
{\cal L}g(\bm{x}) &=&
\frac{1}{2}\sum_{i,j=1}^dx_i(\delta_{ij}-x_j)
\mathbb{E}\Bigl [g_{ij}\big (\bm{x}(1-W)+WV\bm{e}_i\big )\Bigr ]
\nonumber \\
&&~~
+ \frac{\theta}{2}\sum_{i=1}^d\Big (\sum_{j=1}^dp_{ji}x_j - x_i\Big )g_i(\bm{x})
\label{L:LWF:mutation}
\end{eqnarray}
acting on functions in ${\cal C}^2(\Delta)$
maps $d$-dimensional polynomials into polynomials of the same degree, so the right eigenvectors  $\{P_{\bm{n}}(\bm{x})\}$ with eigenvalues $-\lambda_{\bm{n}}$ are polynomials of the same degree satisfying
\begin{equation}
{\cal L}P_{\bm{n}}(\bm{x}) = -\lambda_{\bm{n}}P_{\bm{n}}(\bm{x}).
\label{L:ev}
\end{equation}
The index $\bm{n}$ is $d-1$ dimensional because of the constraint that $\sum_1^dx_j=1$. The eigenvalues $\lambda_{\bm{n}}$, (\ref{L:lambda}) in the following theorem, have a linear form in the $d-1$ non-unit eigenvalues of $I-P$ with coefficients $n_1,\ldots ,n_{d-1}$ which defines $\bm{n}$.
\medskip

\noindent
{\bf Theorem 5.}
\emph{
Let $\{\lambda_{\bm{n}}\}$, $\{P_{\bm{n}}(\bm{x})\}$ be the eigenvalues and right eigenvectors of ${\cal L}$, (\ref{L:LWF:mutation}), satisfying (\ref{L:ev}). Denote the $d-1$ eigenvalues of $P$ which have modulus less than 1 by $\{\phi_k\}_{k=1}^{d-1}$ corresponding to eigenvectors which are rows of a $d-1\times d$ matrix $R$ satisfying 
\begin{equation*}
\sum_{i=1}^dr_{ki}p_{ji} = \phi_kr_{kj},\>k=1,\ldots ,d-1.
%\label{P:ev}
\end{equation*}
Define a $d-1$ dimensional vector $\bm{\xi} = R\bm{x}$. Then the polynomials $P_{\bm{n}}(\bm{x})$ are polynomials in $\bm{\xi}$ whose only leading term of degree $n$ is $\bm{\xi}^{\bm{n}}$ and
\begin{equation}
\lambda_{\bm{n}} = \frac{1}{2}n(n-1)\mathbb{E}\Bigl [(1-W)^{n-2}\Bigr ]
+ \frac{\theta}{2}\sum_{k=1}^{d-1}(1-\phi_k)n_k.
\label{L:lambda}
\end{equation}
}
\medskip

\noindent
{\bf Proof.} 
The second order derivative term in ${\cal L}$ acting on $\bm{x}^{\bm m}$ is
\begin{eqnarray*}
&&-\frac{1}{2}\sum_{i,j=1}^dx_ix_j\mathbb{E}\Bigl [(1-W)^{m-2}\Bigr ]
m_i(m_j-\delta_{ij})\bm{x}^{\bm{m} - \bm{e}_i - \bm{e}_j}
+ \text{~lower~order~terms}
\nonumber \\
&&~~~~~= -\frac{1}{2}m(m-1)\mathbb{E}\Bigl [(1-W)^{m-2}\Bigr ]\bm{x}^{\bm{m}}+\text{~lower~order~terms}.
\end{eqnarray*}
Therefore the same term acting on $\bm{\xi}^{\bm{n}}$ with $m=n$ is
\begin{equation}
-\frac{1}{2}n(n-1)\mathbb{E}\Bigl [(1-W)^{n-2}\Bigr ]\bm{\xi}^{\bm{n}}
+ \text{~lower~order~terms~in~}\bm{\xi}.
\label{lambda:a}
\end{equation}
The linear differential term acting on $\bm{\xi}^{\bm{n}}$ is 
\begin{eqnarray}
&&\frac{\theta}{2}\sum_{i=1}^d\Big (\sum_{j=1}^dp_{ji}x_j - x_i\Big )\frac{\partial}{\partial x_i}
\bm{\xi}^{\bm{n}}
\nonumber \\
&=&
\frac{\theta}{2}\sum_{k=1}^{d-1}
\sum_{i=1}^d\Big (\sum_{j=1}^dp_{ji}x_j - x_i\Big )n_kr_{ki}
\bm{\xi}^{\bm{n}-\bm{e}_k}
\nonumber \\
&=&
-\frac{\theta}{2}\sum_{k=1}^{d-1}(1-\phi_k)n_k\bm{\xi}^{\bm{n}}.
\label{lambda:b}
\end{eqnarray}
(\ref{lambda:a}) and (\ref{lambda:b}) are enough to complete the proof of (\ref{L:lambda}). Suppose we have constructed $\{P_{\bm{m}}(\bm{x})\}_{m < n}$. Then take
\[
P_{\bm{n}}(\bm{x}) = \bm{\xi}^{\bm{n}} - \sum_{\bm{m}: m < n}a_{\bm{n}\bm{m}}P_{\bm{n}}(\bm{x})
\]
where the coefficients are to be determined.
\[
{\cal L}P_{\bm{n}}(\bm{x}) = -\lambda_{\bm{n}}P_{\bm{n}}(\bm{x}) + \sum_{\bm{m}: m < n}b_{\bm{n}\bm{m}}P_{\bm{m}}(\bm{x}) - \sum_{\bm{m}: m < n}a_{\bm{n}\bm{m}}\lambda_{\bm{m}}P_{\bm{m}}(\bm{x})
\]
for determined constants $b_{\bm{n}\bm{m}}$. Choosing 
$
a_{\bm{n}\bm{m}}\lambda_{\bm{m}} = b_{\bm{n}\bm{m}}
$
completes the construction.
\qed
\medskip

\noindent
{\bf Corollary 3.}
\emph{The generator (\ref{L:WFM}) with no mutation term has eigenvalues
\begin{equation*}
\lambda_{\bm{n}} = \lambda_{n} = \frac{1}{2}n(n-1)\mathbb{E}\Bigl [(1-W)^{n-2}\Bigr ]
%\label{Cev:a}
\end{equation*}
repeated ${n+d-2\choose n}$ times and eigenfunctions $\{P_{\bm{n}}(\bm{x})\}_{n\geq 2}$.
}
\qed
\medskip

\noindent
{\bf Corollary 4.}
\emph{
In the parent independent model of mutation the generator has eigenvalues
\begin{equation}
\lambda_{\bm{n}} = \lambda_{n} = \frac{1}{2}n\left \{(n-1)\mathbb{E}\Bigl [(1-W)^{n-2}\Bigr ]
+ \theta\right \}
\label{Cev:b}
\end{equation}
repeated ${n+d-2\choose n}$ times and eigenfunctions $\{P_{\bm{n}}(\bm{x})\}_{n\geq 1}$.
}
\medskip

\noindent
{\bf Proof.} The transition matrix $P$ has rows $(\theta_1/\theta,\ldots ,\theta_d/\theta)$. The right eigenvectors of $P$ are one vector of units with eigenvalue 1, and $d-1$ other vectors such that $\sum_{i=1}^dr_{ki}\theta_i/\theta = 0$. Thus $\phi_k=0$, $k=1,\ldots ,d-1$ and $\lambda_{\bm{n}}$ is equal to (\ref{Cev:b}).
\qed
\medskip

In two dimensions the generator is specified by
\begin{equation}
{\cal L}g(x) = \frac{1}{2}x(1-x)\mathbb{E}\Bigl [g^{\prime\prime}\big (x(1-W) + WV\big )\Bigr ]+
\frac{1}{2}(\theta_1-\theta x)g^{\prime}(x).
\label{L:WFM:mut}
\end{equation}
The eigenvalues are
\[
\lambda_n = \frac{1}{2}n\left \{(n-1)\mathbb{E}\Bigl [(1-W)^{n-2}\Bigr ] + \theta\right \}
\]
and the eigenvectors are polynomials satisfying
\[
{\cal L}P_n(x) = -\lambda_nP_n(x), \>n \geq 1.
\] 
The eigenvalues and polynomials do not depend on $W,V$ for $n=1,2$.
Writing the eigenvalue equation as
\begin{eqnarray*}
&&x(1-x)\mathbb{E}\big [P_n^{\prime\prime}\big (x(1-W)+VW\big )\big ] + (\theta_1-\theta x)P_n^\prime(x)
\nonumber \\
 &&~~~~~~~~~~~~~+ n\big \{(n-1)\mathbb{E}\big [(1-W)^{n-2}\big ]+\theta\big \}P_n(x) = 0
\end{eqnarray*}
there is a similarity to the hypergeometric equation for the Jacobi polynomials  which are the eigenvectors when $W\equiv 0$ \citep{K1964}. Writing the $n^{\text{th}}$ Jacobi polynomial with index parameters of $(\theta_1,\theta_2)$, orthogonal on the Beta distribution with the same parameters as
$\widetilde{P}_n^{(\theta_1,\theta_2)}(x)\equiv z$ for ease of notation, the hypergeometric equation is
\begin{equation}
x(1-x)z^{\prime\prime} + (\theta_1-\theta x)z^\prime 
 + n\big ((n-1)+\theta\big )z = 0,
\end{equation}
see for example \citet{I2005}. Usually the Jacobi polynomials $P_n^{(\alpha_1,\alpha_2)}(x)$ are defined as orthogonal on the weight function
\[
(1-x)^{\alpha_1}(1+x)^{\alpha_2},\>\> -1 < x < 1.
\]
so the translation to orthogonal polynomials on the Beta $(\theta_1,\theta_2)$ distribution is that
\[
\widetilde{P}_n^{(\theta_1,\theta_2)}(x) = 
P_n^{(\theta_2-1,\theta_1-1)}(2x-1).
\]
\subsection{Stationary distributions}
If the mutation matrix $P$ is recurrent, then there is a stationary distribution for the process with generator (\ref{L:LWF:mutation}).
The first and second order moments do not depend on $W$ because they can be found from the generator equations
\[
\mathbb{E}\big [{\cal L}X_i\Big ]=0,\>\mathbb{E}\big [{\cal L}X_iX_j\Big ]=0
\]
which do not depend on $W$ as the second derivatives of $X_i$ and $X_iX_j$ are constant.

 In particular, for the parent independent model of mutation,
comparing moments with those of the Dirichlet $(\bm{\theta})$ distribution which is the stationary distribution for the Wright-Fisher diffusion we have that for $i,j =0,1,\ldots, d$ for any $F$
\[
\mathbb{E}\Bigl[X_i\Bigr ] = \frac{\theta_i}{\theta}
\text{~and~}
\mathbb{E}\Bigl[X_iX_j\Bigr ] = \frac{\theta_i(\theta_j+\delta_{ij})}{\theta(\theta+1)},
\]
with expectation in the stationary distribution (see for example
\citet{E1972}).
Now consider the simplest case, the stationary distribution in two dimensions when the generator is (\ref{L:WFM:mut}). An interesting recurrence for the moments of $X$, the frequency of the first allele, is found in terms of size-biassed versions of $X$.
\medskip

\noindent
{\bf Theorem 6.}
\emph{
Let $Z$ be a random variable with the size-biassed distribution of $X$, $Z_*$ a size-biassed $Z$ random variable; $Z^*$ a size-biassed random variable with respect to $1-Z$; $W=UY$, where $Y$ has distribution $F$, $U$ has a density $2u$, $u \in (0,1)$; $V$  a uniform random variable on $(0,1)$; $B$ a Bernoulli random variable
such that $P(B=1) = \theta_2/\Big(\theta(\theta_1+1)\Big )$ with $U,V,Y,Z^*,Z_*,B$ independent.
Then
\begin{equation}
VZ_* =^{\cal D} \big (1-B \big)VZ + B\big (Z^*(1-W) + WV\big ).
\label{2dim:0}
\end{equation}
}
\medskip

\noindent{\bf Proof.}
Let $g(x)=x^{n+2}$, then since $\mathbb{E}\Bigl [{\cal L}g(X)\Bigr ] = 0$ with expectation in the stationary distribution 
\begin{equation*}
\frac{(n+2)(n+1)}{2}\mathbb{E}\Bigl [\Bigl (X(1-W) + WV\Bigr )^nX(1-X)\Bigr ]
+\frac{n+2}{2}\mathbb{E}\Bigl [X^nX(\theta_1-\theta X)\Bigr ] = 0
\end{equation*}
or
\begin{equation}
\frac{\theta}{n+1}\mathbb{E}\Bigr [X^nX^2\Bigr ] = \frac{\theta_1}{n+1}\mathbb{E}\Bigr [X^nX\Bigr ]
+\mathbb{E}\Bigl [\Bigl (X(1-W) + WV\Bigr )^nX(1-X)\Bigr ].
\label{moment:eq:0}
\end{equation}
Let $Z$ be a random variable with the size-biassed distribution of $X$, $Z_*$ a size-biassed $Z$ random variable and $Z^*$ a size-biassed random variable with respect to $1-Z$. The distribution of $Z$ is re-weighted by $Z$ and divided by $\mathbb{E}\big [Z\big ]$ to obtain the distribution of $Z_*$ and similarly the distribution is weighted by $1-Z$ and divided by $\mathbb{E}\big [1-Z\big ]$ to obtain the distribution of 
$Z^*$.
Then knowing that 
\[
\mathbb{E}\Bigl [X\Bigr ]=\frac{\theta_1}{\theta},\>
\mathbb{E}\Bigl [X^2\Bigr ]=\frac{\theta_1(\theta_1+1)}{\theta(\theta+1)},\>
\mathbb{E}\Bigl [X(1-X)\Bigr ]=\frac{\theta_1(\theta-\theta_1)}{\theta(\theta+1)}
\]
(\ref{moment:eq:0}) can be written as
\begin{equation}
\mathbb{E}\Bigl [(VZ_*)^n\Bigr ] = \frac{\theta_1(\theta+1)}{\theta(\theta_1+1)}\mathbb{E}\Bigl [(VZ)^n\Bigr ]
+ \frac{\theta - \theta_1}{\theta(\theta_1+1)}\mathbb{E}\Bigl [\Bigl (Z^*(1-W) + WV\Bigr )^n\Bigr ].
\label{moment:eq:1}
\end{equation}
Recall that
\[
P(B=1) = \frac{\theta_2}{\theta(\theta_1+1)}.
\]
Then (\ref{moment:eq:1}) implies the distributional identity (\ref{2dim:0}).
\qed
\medskip

This equation may be related to a renewal process, because the distribution of excess life $\gamma_t$ in a renewal process with increments distributed as $Z$ satisfies
\[
\lim_{t\to \infty}P(\gamma_t > \eta) = P(VZ_* > \eta) = \int_\eta^1\frac{P(Z>z)}{\mathbb{E}[Z]}dz,
\]
where $\mathbb{E}\big [Z\big ] = \theta_1/\theta$ \citep{C1970}.

The identity (\ref{2dim:0}) implies an integral equation for the stationary distribution in the 2-dimensional model.
\medskip

\noindent
{\bf Theorem 7.}
\emph{
Let $f_X(u)$, $0 < u < 1$ be the stationary density in the diffusion process with generator (\ref{L:WFM}), and $f_W(w)$ be the density of $W$. Suppose that $F$ has no atom at zero. Then $f_X(u)$ satisfies the integral equations
\begin{eqnarray}
\big (\theta_1 - \theta u\big )f_X(u)
&=&
-\int_0^u\frac{1}{u-z}f_W\Big (1-\frac{1-u}{1-z}\Big )z(1-z)f_X(z)dz
\nonumber \\
&&~~+\int_u^1\frac{1}{z-u}f_W\Big (1-\frac{u}{z}\Big )z(1-z)f_X(z)dz
\label{statThm:0}
\end{eqnarray}
and
\begin{eqnarray}
\big (\theta_2 - \theta (1-u)\big )f_X(u)
&=&
\int_0^u\frac{1}{u-z}f_W\Big (1-\frac{1-u}{1-z}\Big )z(1-z)f_X(z)dz
\nonumber \\
&&~~-\int_u^1\frac{1}{z-u}f_W\Big (1-\frac{u}{z}\Big )z(1-z)f_X(z)dz.
\label{statThm:1}
\end{eqnarray}
These equations are equivalent to
\begin{eqnarray}
\big (\theta_1 - \theta u\big )f_X(u)
&=&
-\int_0^u2F^+\Big (1-\frac{1-u}{1-z}\Big )zf_X(z)dz
\nonumber \\
&&~+\int_u^12F^+\Big (1-\frac{u}{z}\Big )(1-z)f_X(z)dz
\label{statThm:0a}
\end{eqnarray}
and
\begin{eqnarray}
\big (\theta_2 - \theta (1-u)\big )f_X(u)
&=&
\int_0^u2F^+\Big (1-\frac{1-u}{1-z}\Big )zf_X(z)dz
\nonumber \\
&&~-\int_u^12F^+\Big (1-\frac{u}{z}\Big )(1-z)f_X(z)dz.
\label{statThm:1a}
\end{eqnarray}
}
\medskip

\noindent
{\bf Proof.}
Let the random line $L=Z^*(1-W) + WV$ as a function of $W$. The line segment $L$ varies from $\min(Z^*,V)$ to $\max(Z^*,V)$ as $W$ varies. The density of the line $L$ conditional on $(Z^*,V) = (z,v)$ is, for $\min(z,v) < u < \max(z,v)$,
\begin{equation*}
f_{L\mid (z,v)}(u) = \frac{1}{|z-v|}f_W\Bigl (\frac{z-u}{z-v}\Bigr )
\end{equation*}
and there is a possible atom
\[
P(L=z\mid (z,v)\big ) = P(W=0).
\]
Splitting the region by $v < z$ and $v > z$, the unconditional density of $L$ is
\begin{eqnarray*}
f_L(u)&=&
P(W=0)f_{Z^*}(u)
\nonumber \\
&&+\int_{0 < v < u < z <1}\frac{1}{z-v}f_W\Big (\frac{z-u}{z-v}\Bigr )f_{Z^*}(z)dzdv
\nonumber \\
&&+\int_{0 < z < u < v <1}\frac{1}{v-z}f_W\Big (\frac{u-z}{v-z}\Bigr )f_{Z^*}(z)dzdv.
\end{eqnarray*}
Changing variables in the integral
\begin{eqnarray}
f_L(u)&=&
P(W=0)f_{Z^*}(u)
\nonumber \\
&&+\int_0^u\int_{1-\frac{1-u}{1-z}}^1\frac{1}{\xi}f_W(\xi)d\xi\> f_{Z^*}(z)dz
\nonumber \\
&&+\int_u^1\int_{1 - \frac{u}{z}}^1\frac{1}{\xi}f_W(\xi)d\xi\> f_{Z^*}(z)dz.
\label{L:density}
\end{eqnarray}
The density identity equivalent to the identity (\ref{2dim:0}) is therefore
\begin{equation}
f_{VZ_*}(u) = P(B=0)f_{VZ}(u) + P(B=1)f_L(u).
\label{density:id:0}
\end{equation}
Note that if $\zeta$ is a random variable on $(0,1)$ with density $f_\zeta(y)$ then the density of $V\zeta$, where $V$ is independent of $\zeta$ and uniform on $(0,1)$ is
\[
f_{V\zeta}(u) = \int_u^1y^{-1}f_\zeta(y)dy.
\]
Therefore (\ref{density:id:0}) is equivalent to
\begin{equation}
\frac{\int_u^1yf_X(y)dy}{\mathbb{E}\bigl [X^2\bigr ]}
= P(B=0)\frac{\int_u^1f_X(y)dy}{\mathbb{E}\bigl [X\bigr ]} + P(B=1)f_L(u).
\label{density:id:1}
\end{equation}
Differentiating (\ref{density:id:1}) the density $f_X(u)$ satisfies the integral equation
\begin{equation}
uf_X(u) = \frac{\theta_1}{\theta}f_X(u) - \frac{1}{\theta}f_{\diamond L}^\prime(u),
\label{integral:eq:0}
\end{equation}
where 
\[
f_{\diamond L}(u)=\mathbb{E}\Bigl [X(1-X)\Bigr ]f_L(u).
\]
$f_{\diamond L}(u)$ is similar to (\ref{L:density}) with $f_{Z_*}(z)$ replaced by $z(1-z)f_X(z)$.
A straightforward calculation gives that when $W$ has no atom at zero
\begin{eqnarray}
f^\prime_{\diamond L}(u) &=& 
-\int_0^u\frac{1}{u-z}f_W\Big (1-\frac{1-u}{1-z}\Big )z(1-z)f_X(z)dz
\nonumber \\
&&~~~~+\int_u^1\frac{1}{z-u}f_W\Big (1-\frac{u}{z}\Big )z(1-z)f_X(z)dz.
\label{L:diff}
\end{eqnarray}
Recalling (\ref{WY:1}) another form is
\begin{eqnarray}
f^\prime_{\diamond L}(u) &=& 
-\int_0^u2F^+\Big (1-\frac{1-u}{1-z}\Big )zf_X(z)dz
\nonumber \\
&&~~~~+\int_u^12F^+\Big (1-\frac{u}{z}\Big )(1-z)f_X(z)dz.
\label{L:diff:a}
\end{eqnarray}
Considering $1-X$ a second integral equation is
\begin{equation}
(1-u)f_X(u) = \frac{\theta_2}{\theta}f_X(u) + \frac{1}{\theta}f_{\diamond L}^\prime(u).
\label{integral:eq:1}
\end{equation}
Substituting the expression (\ref{L:diff}) for $f^\prime_{\diamond L}(u)$ in (\ref{integral:eq:0}) and (\ref{integral:eq:1}) gives (\ref{statThm:0}) and (\ref{statThm:1}). The alternative form (\ref{L:diff:a}) gives (\ref{statThm:0a}) and (\ref{statThm:1a}).
\qed
\medskip

Another approach that imitates the usual way of finding the stationary distribution in a diffusion process is to consider the equation
\begin{equation}
\int_0^1\big [{\cal L}g(x)\big ]f_X(x)dx = 0,
\label{stationary:0}
\end{equation}
where $g(x)$ is a test function in ${\cal C}^2([0,1])$. Denote $\sigma^2(x) = x(1-x)$, $\mu(x) = \theta_1 - \theta x -\sigma x(1-x)$ and let
\[
k(x) = \mathbb{E}\Bigl [(1-W)^{-2}g\big (x(1-W) + VW\big )\Bigr ].
\]
Equation (\ref{stationary:0}) is equivalent to 
\begin{equation}
\int_0^1\Bigl [\frac{1}{2}\sigma^2(x)\frac{d^2}{dx^2}k(x) + \mu(x)\frac{d}{dx}g(x)\Bigr ]f_X(x)dx = 0.
\label{stationary:1}
\end{equation}
Integrating by parts in (\ref{stationary:1}) and taking care with boundary conditions gives that
\begin{eqnarray*}
&&0=\int_0^1\Bigl [k(x)\frac{1}{2}\frac{d^2}{dx^2}\bigl [\sigma^2(x)f_X(x)\bigr ] - g(x)\frac{d}{dx}\bigl [\mu(x)f_X(x)\bigr ]\Bigr ]dx
\nonumber\\
&&+\Biggl [\big [\frac{d}{dx}k(x)\big ]\big [\frac{1}{2}\sigma^2(x)f_X(x)\big ] -k(x)\frac{d}{dx}\big [\frac{1}{2}\sigma^2(x)f_X(x)\bigr ]
+ g(x)\mu(x)f_X(x)\Biggr ]_0^1.
\end{eqnarray*}
If $W\equiv 0$ then $k(x) = g(x)$ and we can conclude that $f_X(x)$ satisfies the forward equation,
\[
\frac{1}{2}\frac{d^2}{dx^2}\bigl [\sigma^2(x)f_X(x)\bigr ] - \frac{d}{dx}\bigl [\mu(x)f_X (x)\bigr ] = 0.
\]
An equivalent approach seems difficult when $k(x)\ne g(x)$.
\subsection{Green's function}
The Green's function $G(x,\xi)$, whether there is mutation and selection or not, is obtained in a standard approach by solving, for $\gamma(x)$, the differential equation
\begin{equation}
{\cal L}\gamma(x) = -g(x),\>\>\gamma(0)=\gamma(1)=0.
\label{Green:1}
\end{equation}
Then
\begin{equation*}
\gamma(x) = \int_0^1G(x,\xi)g(\xi)d\xi.
%\label{Green:2}
\end{equation*}
Consider the model with no selection.
Equation (\ref{Green:1}) is non-linear, equivalent to
\begin{equation*}
\frac{1}{2}x(1-x)\mathbb{E}\Bigl [\gamma^{\prime\prime}(x(1-W)+VW)\Bigr ] + \frac{1}{2}(\theta_1-\theta x)\gamma^{\prime}(x) = -g(x),
%\label{Green:3}
\end{equation*}
or
\begin{equation}
\frac{1}{2}x(1-x)k^{\prime\prime}(x) + \frac{1}{2}(\theta_1-\theta x)\gamma^{\prime}(x) = -g(x),
\label{Green:4}
\end{equation}
where
\begin{equation*}
k(x) = \mathbb{E}\Bigl [(1-W)^{-2}\gamma\big (x(1-W) + VW\big )\Bigr ].
\end{equation*}
In the simplest case when $\theta=0$, (\ref{Green:4}) becomes
\begin{equation*}
k^{\prime\prime}(x) = -2\frac{g(x)}{x(1-x)}.
\end{equation*}
Taking a standard Green's function approach, with care that $k(0),k(1)$ are not zero,
\begin{eqnarray*}
k(x) &=& k(0)(1-x)  + k(1)x
\nonumber \\
&&~~ + (1-x)\int_0^x\frac{2g(\eta)}{1-\eta}d\eta
+ x\int_x^1\frac{2g(\eta)}{\eta}d\eta.
\end{eqnarray*}
If $g(x)=1,\>x \in (0,1)$ then $\gamma (x)$ is the mean time to absorption at 0 or 1 when $X(0)=x$. There is a non-linear equation to solve of
\begin{equation}
k(x) = k(0)(1-x) + k(1)x - 2(1-x)\log (1-x) - 2 x\log x.
\label{nonLinear:k}
\end{equation}
It is possible that $\gamma(x)=\infty$ if the $\Lambda$-coalescent does not come down from infinity.
\subsection{The frequency spectrum in the infinitely-many-alleles model}
We consider the infinitely-many-alleles model as a limit from a $d$-allele model with $\theta_i=\theta/d$, $i=1,\ldots,d$. The limit is thought of as a limit from $d$ points $X^d_1,\ldots ,X^d_d~$ to points of a point process $\{X_i\}_{i=1}^\infty$. The 1-dimensional frequency spectrum $\mu$ is a non-negative measure such that for functions $f$ in ${\cal C}([0,1])$ such that $f(x)/x$ is bounded as $x \to 0$, with expectation in the stationary distribution
\begin{equation}
\mathbb{E}\Bigl [\sum_{i=1}^\infty f(X_i)\Bigr ] = \int_0^1f(x)\mu(dx).
\label{fsdef:0}
\end{equation}
There is an assumption that the point process does not have multiple points at any single position for (\ref{fsdef:0}) to hold.
Symmetry in the $d$-allele model shows that
\[
\int_0^1f(x)\mu(dx) = \lim_{d\to \infty}d\mathbb{E}\Bigl [f(X_1)\Bigr ].
\]
\textcolor{black}{If the $\Lambda$-coalescent does not come down from infinity then there may be an accumulation of points at zero as $d \to \infty$ and $\int_0^1 x\mu(x)dx < 1$.  We do not consider this case in the next theorem.}
The classical Wright-Fisher diffusion gives rise to the Poisson-Dirichlet process with a frequency spectrum of
\[
\mu(dx)=\theta x^{-1}(1-x)^{\theta - 1}dx,\> 0 < x < 1.
\] 
\medskip

\noindent
{\bf Theorem 8.}
\emph{Let $\mu(dz)$ be the frequency spectrum measure \textcolor{black}{in an infinitely-many-alleles $\Lambda$-Fleming-Viot process which comes down from infinity} and $Z$ a random variable with probability measure $z\mu(dz)$. Let $Z_*$ be a random variable with a size-biassed distribution of $Z$ and $Z^*$ a random variable with a size-biassed distribution of $Z$ with respect to $1-Z$. $Z_*$ has a measure $(\theta + 1)z^2\mu(dz)$ and $Z^*$ has a measure $\theta^{-1}(\theta+1)z(1-z)\mu(dz)$, $0 < z < 1$. Then
\begin{equation}
VZ_* =^{\cal D} Z^*(1-W) + WV,
\label{freq:id:0}
\end{equation}
where $V,Z_*,Z^*, W$ are independent. Let $\mu(dz) = \beta (z)dz$. Suppose that $Y$ has no atom at zero. 
Then an integral equation for $\beta(x)$ is
\begin{eqnarray*}
\theta u\beta (u) &=&
\int_0^u\frac{1}{u-z}f_W\Big (1-\frac{1-u}{1-z}\Big )z(1-z)\beta (z)dz
\nonumber \\
&&~~~~-\int_u^1\frac{1}{z-u}f_W\Big (1-\frac{u}{z}\Big )z(1-z)\beta (z)dz
\label{fsinteg:0}
\end{eqnarray*}
which is equivalent to
\begin{eqnarray*}
\theta u\beta (u) &=&
 \int_0^u2F^+\Big (1-\frac{1-u}{1-z}\Big )z\beta (z)dz
\nonumber \\
&&~~~~-\int_u^12F^+\Big (1-\frac{u}{z}\Big )(1-z)\beta (z)dz.
\label{fsinteg:0a}
\end{eqnarray*}
}
\medskip

\noindent
{\bf Proof.}
To obtain a limit in the $\Lambda$-Fleming-Viot process let $\theta_1=\theta/d$ in the generator (\ref{L:WFM:mut}). In the identity (\ref{2dim:0}) the density of $Z$ is $d\, zf_{X_1}(z)$, $0 < z < 1$, by symmetry. Let $d\to \infty$ in the identity (\ref{2dim:0}). Then the identity becomes (\ref{freq:id:0}).
%The constant $\theta$ appears in the identity (\ref{freq:id:0}) because of scaling in the size-biassed distributions. $Z_*$ has a measure $(\theta + 1)z^2\mu(dz)$ and $Z^*$ has a measure $\theta^{-1}(\theta+1)z(1-z)\mu(dz)$, $0 < z < 1$.
%Let $\mu(dz) = \beta (z)dz$.
%%
The integral equations for the stationary distribution when there are two types imply an integral equation for $\beta (x)$. In view of (\ref{integral:eq:0}) 
\begin{equation}
x\beta (x) = -\frac{1}{\theta}f_{\bullet L}^\prime(u),
\end{equation}
where $f_{\bullet L}^\prime (u)$ is similar to $f_{\diamond L}^\prime (u)$ with $f_X(z)$ replaced by $\beta(z)$.
\qed
\medskip

In the Wright-Fisher diffusion $W\equiv 0$ and the identity (\ref{freq:id:0}) is
\begin{equation}
VZ_* =^{\cal D}Z^*.
\label{freq:id:1}
\end{equation}
It is straightforward to verify that if $Z$ has density $\theta (1-z)^{\theta-1}$, then (\ref{freq:id:0}) is satisfied. A direct solution can be found in the following way. From  (\ref{freq:id:1}) 
\begin{equation}
\theta \int_z^1y\beta(y)dy = z(1-z)\beta(z),
\label{basic}
\end{equation}
where $\theta$ is defined by
\[
\theta = \frac{ \int_0^1z^2\beta (z)dz}{ \int_0^1z(1-z)\beta (z)dz}.
\]
Write (\ref{basic}) as
\[
\frac{d}{dz}\log \int_z^1y\beta(y)dy = -\theta (1-z)^{-1}.
\]
Solving this differential equation,
\[
\log \int_z^1y\beta (y)dy = \theta \log (1-z) + A,
\]
for a constant $A$. Therefore
\[
\int_z^1y\beta (y)dy = (1-z)^\theta
\]
because $\int_0^1yh(y)dy = 1$, and since (\ref{basic}) holds
\[
\beta(z) = \theta z^{-1} (1-z)^{\theta-1},\> 0 < z < 1.
\]
\subsection{A different dual process}
The typed $\Lambda$-coalescent tree process is a moment dual in the Fleming-Viot process, see for example \citet{EGT2010}.
We work through a different type of dual process which is a death process decreasing in steps of 1. Let  $d=2$ for simplicity. The generator ${\cal L}$ is specified by (\ref{L:WFM:mut}).
Let $\{g_n(x)\}$ be a sequence of
monic polynomials that are defined below satisfying the
generator equation
\begin{eqnarray}
{\cal L}g_n &=&
\frac{1}{2}x(1-x)\mathbb{E}g_n^{\prime\prime}(x(1-W)+VW) +
\frac{1}{2}(\theta_1-\theta x)g_n^{\prime}(x)\nonumber \\
&=& {n\choose 2}\mathbb{E}(1-W)^{n-2}[g_{n-1}(x)-g_n(x)] +
n\frac{1}{2}[\theta_1 g_{n-1}(x) - \theta g_n(x)]
\nonumber \\
\label{dual:0}
\end{eqnarray}
with $g_0(x) = 1$. Equation (\ref{dual:0}) is an analogue of the
Wright-Fisher diffusion when we look at $g_n(x) = x^n$, with the
second line chosen to mimic the Wright-Fisher case. Rearrange the equation to
define $g_n(x)$ in terms of $g_{n-1}(x)$ as
\begin{eqnarray}
&&\frac{1}{2}x(1-x)\mathbb{E}g_n^{\prime\prime}(x(1-W)+VW) +
\frac{1}{2}(\theta_1-\theta x)g_n^{\prime}(x) +\lambda_n g_n(x)\nonumber \\
&&~~= \frac{n}{2}\Bigl [(n-1)\mathbb{E}(1-W)^{n-2} +\theta_1\Bigr ]g_{n-1}(x) 
\label{dual:1}
\end{eqnarray}
The polynomials $\{g_n(x)\}$ are well defined by (\ref{dual:1}) by
recursively calculating the coefficients of $x^r$ in $g_n(x)$ from
$r=n-1,n-2,\ldots ,0$.
\medskip

\noindent
{\bf Theorem 9.}
\emph{Let $\{g_n(x)\}_{n=0}^\infty$ defined by
(\ref{dual:1}). If $X$ has a stationary distribution then
\begin{equation}
\mathbb{E}\Big [g_n(X)\Big ]
= \frac
{\prod_{j=1}^{n}\Big ((j-1)\mathbb{E}(1-W)^{j-2} + \theta_1\Big ) }
{\prod_{j=1}^{n}\Big ((j-1)\mathbb{E}(1-W)^{j-2} + \theta\Big ) }.
\label{Stationaryg:0}
\end{equation} 
Let $h_n(x) = g_n(x)/\mathbb{E}\Big [g_n(X)\Big ]$. There is a dual process $\{N(t)\}_{t\geq 0}$ to $\{X(t)\}_{t \geq 0}$ based on the test functions $\{h_n(x)\}_{n=0}^\infty$ which is a death process with rates $n\to n-1$, $n \geq 1$, of
\[
\lambda_n = \frac{n}{2}\Bigl [(n-1)\mathbb{E}(1-W)^{n-2} + \theta\Bigr ].
\]
The dual equation is
\begin{equation}
\mathbb{E}_{X(0)=x}\Bigl [h_n(X(t))\Bigr ] = \mathbb{E}_{N(0)=n}\Bigl [h_{N(t)}(x)\Bigr ].
\label{dual:5}
\end{equation}
The transition functions for the dual process are
\begin{eqnarray}
&&P\bigl (N(t) = j\mid N(0)=i \big )
\nonumber \\
&&~~~= \sum_{k=j}^ie^{-\lambda_k t}r_i^{(k)}l_j^{(k)}
\nonumber \\
&&~~~=
\sum_{k=j}^ie^{-\lambda_k t}(-1)^{k-j}\frac{\lambda_{j+1}\cdots \lambda_i}
{(\lambda_j - \lambda_k)\cdots (\lambda_{k-1}-\lambda_k)(\lambda_{k+1}-\lambda_k)\cdots (\lambda_i - \lambda_k)}.
\nonumber \\
\label{expansion:0}
\end{eqnarray}
$\{N(t)\}_{t\geq 0}$ comes down from infinity if and only if
\begin{equation}
\int_1^\infty\frac{dq}{q^2\mathbb{E}\big [(1-W)^q\big ]} < \infty.
\label{come:a0}
\end{equation}
which implies the $\Lambda$-coalescent coming down from infinity.
The distribution of $N(t)$ given an entrance boundary at infinity is \begin{equation}
P\bigl (N(t) = j\mid N(0)=\infty \big ) = \sum_{k=j}^\infty e^{-\lambda_k t}r_\infty^{(k)}l_j^{(k)},
\label{expansion:1}
\end{equation}
where 
\[
r_\infty^{(k)} = \prod_{l=k}^\infty \Bigl (1 - \frac{\lambda_k}{\lambda_l}\Bigr )^{-1},
\]
well defined assuming the condition (\ref{come:a0}) when the coalescent comes down from infinity.
}
\medskip

\noindent
{\bf Proof.}
Equation (\ref{Stationaryg:0}) follows directly from 
$\mathbb{E}\big [g_n(X)\big ]$ in (\ref{dual:0}).
Note that
\begin{equation*}
{\cal L}h_n = \lambda_n[h_{n-1}-h_n].
\label{dual:4}
\end{equation*}
which is correctly set up as a dual generator equation of the death process $\{N(t)\}_{t \geq 0}$. The dual equation is then (\ref{dual:5}).

The process $\{N(t),t\geq 0\}$ comes down from infinity if and only if
\begin{equation}
\sum_2^\infty \lambda_n^{-1} < \infty,
\label{come:a}
\end{equation}
which implies the $\Lambda$-coalescent coming down from infinity because
\[
\lambda_n = \sum_{k=2}^n{n\choose k}\lambda_{nk} + n\theta
\]
so (\ref{come:a}) is equivalent to
\begin{equation}
\sum_{n=2}^\infty \Biggl[\sum_{k=2}^n{n\choose k}\lambda_{nk}\Biggr ]^{-1} < \infty,
\label{come:b}
\end{equation}
and
\begin{equation}
\sum_{n=2}^\infty \Biggl [\sum_{k=2}^n(k-1){n\choose k}\lambda_{nk}\Biggr ]^{-1} < 
\>\>\sum_{n=2}^\infty \Biggl[\sum_{k=2}^n{n\choose k}\lambda_{nk}\Biggr ]^{-1} < \infty.
\label{series:00}
\end{equation}
Recalling that 
\[
\sum_{k=2}^n{n\choose k}\lambda_{nk} 
= \frac{1}{2}n(n-1)\mathbb{E}\big [(1-W)^{n-2} \big ],
\]
by the integral comparison test (\ref{come:b}) is equivalent to 
(\ref{come:a0}).
\qed
\medskip

For example if $W$ has a Beta $(\alpha,\beta)$ distribution for $\alpha,\beta >0$ then 
$
\mathbb{E}\bigl [(1-W)^n\bigr ] \sim Cn^{-\alpha},
$
where $C$ is a constant, so if $\alpha < 1$, then
$
\sum_{n=2}^\infty \lambda_n^{-1} < \infty,
$
because the $n$th term is asymptotic to 
$
(C/2)n^{2-\alpha}
$. Coming down from infinity does not necessarily imply that (\ref{come:a}) or (\ref{come:b}) hold. 
%Example 14 from \citet{S2000} where the $\Lambda$ measure has a density $-\log (y)$, $y\in (0,1)$ is such a case. This coalescent comes down from infinity, however $\lambda_n \sim n\log(n)$ so (\ref{come:b}) does not hold.
%
In general the tail of the series (\ref{series:00}) 
\[ 
\sum_{n=N}^\infty\Biggl[\sum_{k=2}^n{n\choose k}\lambda_{nk}\Biggr ]^{-1}
\approx 
\frac{1}{2}\int_N^\infty \frac{1}{q^2}
\frac{dq}
{\mathbb{E}\big [(1-W)^q\big ]}
= \frac{1}{2}\int_0^{N^{-1}}\frac{dz}
{\mathbb{E}\big [(1-W)^{z^{-1}}\big ]}.
\]
Convergence of the integral depends on $\mathbb{E}\big [(1-W)^{z^{-1}}\big ]$ being large enough as $z \to 0$. 
It is very likely that there are connections with the speed of coming down from infinity studied in \citet{BBLa2012}, but the exact connections are not clear.

The transition functions for the process $\{N(t),t\geq 0\}$ are easily found from an eigenfunction analysis of the $Q$ matrix, where $q_{jj}=-\lambda_j$ and $q_{jj-1} =  \lambda_j$. The approach in finding the eigenfunction expansion for the transition distribution in the Kingman coalescent is taken in \citet{T1984} (see also \citet{G1980}). 
The left and right eigenvectors $l_j^{(k)}$ and $r_i^{(k)}$ are triangular in form with $l_j^{(k)}=0,\>j>k$
and $r_i^{(k)}=0,\>i < k$. Explicit formulae are $l_k^{(k)}=r_k^{(k)} = 1$ and
\begin{eqnarray}
l_j^{(k)}&=& \frac{(-1)^{k-j}\lambda_{j+1}\cdots \lambda_{k}}
{(\lambda_j-\lambda_k)\cdots (\lambda_{k-1} - \lambda_{k})},\>j < k,
\nonumber \\
r_i^{(k)} &=& \frac{\lambda_i\cdots \lambda_{k+1}}
{(\lambda_i-\lambda_k)\cdots (\lambda_{k+1}-\lambda_k)},\>i > k.
\label{lrev}
\end{eqnarray}
The transition functions are then given by (\ref{expansion:0}).

The distribution of $N(t)$ given an entrance boundary at infinity is the distribution as $i \to \infty$ which is (\ref{expansion:1}),
well defined assuming the condition (\ref{come:a}) when the coalescent comes down from infinity.
\qed
\medskip

The condition of \citet{BLG2006}, (\ref{comedown:0}), for coming down from infinity is equivalent to
\begin{equation}
 \int_1^\infty \frac{dq}{q\mathbb{E}\Big [\frac{1-e^{-qW}}{W}\Big ]}
 < \infty.
\label{imply:0}
\end{equation}
There can be a gap where the $\Lambda$-coalescent comes down from infinity but $\{N(t)\}_{t\geq 0}$ does not come down from infinity because (\ref{come:a0}) and (\ref{imply:0}) are not equivalent.
\subsubsection*{Eigenfunctions $P_n(x)$ and polynomials $g_n(x)$}
It is extremely interesting that 
the polynomials $\{P_n(x)\}$ are analogous to the monic Jacobi polynomials distribution with $\{g_n(x)\}$ analogous to $\{x^n\}$. 

Express
\[
g_n(x) = P_n(x) + \sum_{r=0}^{n-1}b_{nr}P_r(x),
\]
where $P_n(x)$ are the eigenfunctions of ${\cal L}$.
Denote
\[
\lambda_n^\circ = \frac{n}{2}\Bigl [(n-1)\mathbb{E}(1-W)^{n-2} + \theta_1\Bigr ].
\]
From (\ref{dual:1}) and noting that
\begin{equation}
{\cal L}P_n = -\lambda_n P_n,\>\>{\cal L}g_n = -\lambda_ng_n + \lambda_n^\circ g_{n-1}
\label{poly:2}
\end{equation}
it follows that
\begin{equation}
\sum_{r=0}^{n-1}b_{nr}[-\lambda_r + \lambda_n]P_r(x) =
\lambda_n^\circ\sum_{r=0}^{n-1}b_{n-1 r}P_r(x). 
\label{dual:2}
\end{equation}
$g_l(x)$ being a monic polynomial means that
$b_{ll}=1$, $l=0,1,\ldots$. Calculating coefficients from
(\ref{dual:2})
\begin{eqnarray*}
b_{nr} &=& \frac{\lambda_n^\circ}{\lambda_n-\lambda_r}b_{n-1 r}
\nonumber \\
&=& \frac{\lambda_n^\circ\lambda_{n-1}^\circ\cdots \lambda_{r+1}^\circ}
{(\lambda_n-\lambda_r)(\lambda_{n-1}-\lambda_r)\cdots
(\lambda_{r+1}-\lambda_r)}.
%\label{dual:3}
\end{eqnarray*}
The eigenfunctions $\{P_n(x)\}$ also have an expansion in terms of the polynomials $\{g_r(x)\}$. Let
\begin{equation}
P_n(x) = g_n(x) + \sum_{r=0}^{n-1}c_{nr}g_r(x).
\label{inv:0}
\end{equation}
From (\ref{poly:2})
\[
-\lambda_nP_n(x) = -\lambda_ng_n(x) + \lambda_n^\circ g_{n-1}(x) + \sum_{r=0}^{n-1}
\Bigl [-\lambda_rc_{nr} + \lambda_{r+1}^\circ c_{nr+1}\Bigr ]g_r(x).
\]
Expressing the left side by the expansion (\ref{inv:0}) and equating coefficients of $g_r(x)$
\[
-\lambda_nc_{nr} = -\lambda_rc_{nr} + \lambda_{r+1}^\circ c_{nr+1}.
\]
The coefficients therefore are
\begin{equation}
c_{nr} = \frac
{ \lambda_{r+1}^\circ\cdots\lambda_n^\circ }
{ (\lambda_r-\lambda_n)\cdots (\lambda_{n-1}-\lambda_n)}.
\label{coeft:c}
\end{equation}

Scale the equation (\ref{dual:0}) by taking
\begin{eqnarray*}
g_n(x) &=& \frac{\lambda_n^\circ\cdots \lambda_1^\circ}{\lambda_n\cdots
\lambda_1}h_n(x)\nonumber \\
&=& \frac
{\prod_{j=1}^{n}\Big ((j-1)\mathbb{E}(1-W)^{j-2} + \theta_1\Big ) }
{\prod_{j=1}^{n}\Big ((j-1)\mathbb{E}(1-W)^{j-2} + \theta\Big ) }h_n(x).
\end{eqnarray*}
Denote $\omega_n$ as a Beta moment analog
\begin{equation*}
\omega_n =  \frac
{\prod_{j=1}^{n}\Big ((j-1)\mathbb{E}(1-W)^{j-2} + \theta_1\Big ) }
{\prod_{j=1}^{n}\Big ((j-1)\mathbb{E}(1-W)^{j-2} + \theta\Big ) }.
\end{equation*}
so
\[
g_n(x) = \omega_nh_n(x).
\]
Note that if $X$ has a stationary distribution then
\[
\mathbb{E}\Bigl [g_n(X)\Bigr ] = \omega_n.
\]
The polynomials $\{P_n(x)\}$ are analogous to the monic Jacobi polynomials orthogonal on the Beta $(\theta_1,\theta_2)$ distribution with $\{g_n(x)\}$ analogous to $\{x^n\}$. If $W\equiv 0$ then they are identical in the analogy. In the Jacobi polynomial case (\ref{coeft:c}) simplifies to
\begin{equation*}
c_{nr} = (-1)^{n-r}\frac{(n-r-1)!}{r!}\frac{{\theta_1}_{(n)}}{{\theta_1}_{(r)}}\frac{(n+\theta)_{(r)}}{(n+\theta)_{(n-1)}}.
%\label{coeft:J}
\end{equation*}
The process is not reversible, so the polynomials are not orthogonal on any measure unless they are the Jacobi polynomials. 
\subsubsection*{Higher dimensions} Let ${\cal L}$ be the $d$-dimensional $\Lambda$-Fleming-Viot generator with mutation and define polynomials $\{g_{\bm{n}}(\bm{x})\}$ with $g_{\bm{0}}(\bm{x})=1$ by
\begin{equation}
{\cal L}g_{\bm{n}}(\bm{x}) = -\lambda_{\bm{n}}g_{\bm{n}}(\bm{x})
+ \frac{1}{2}\sum_{i=1}^d\frac{n_i}{n}\cdot n\Bigl ((n_i-1)\mathbb{E}\Bigl [(1-W)^{n_i-2}\Bigr ] + \theta_i\Bigr )g_{\bm{n}-\bm{e}_i}(\bm{x}).
\end{equation}
This is an analogy with the Wright-Fisher generator acting on $\bm{x}^{\bm{n}}$.
The polynomials are well defined by recursion on their coefficients. In a similar calculation to the two dimensional case there is a Dirichlet moment analogue
\begin{equation}
\mathbb{E}\big [g_{\bm{n}}(\bm{X}) \big ]
= \frac{
\prod_{i=1}^d\Bigl [\prod_{j=1}^{n_i}\bigl ((j-1)\mathbb{E}\bigl [(1-W)^{j-2}\bigr ]+\theta_i\bigr )\Bigr ]
}
{
\prod_{j=1}^{n}\bigl ((j-1)\mathbb{E}\bigl [(1-W)^{j-2}\bigr ]+\theta\bigr )
}.
\end{equation}
The dual process constructed from test functions $g_{\bm{n}}(\bm{X})/\mathbb{E}\bigl [g_{\bm{n}}(\bm{X})\bigr ]$ has transitions 
\[
\bm{n} \to \bm{n}-\bm{e}_i\text{~at~rate~} \frac{n_i}{n}\Bigl ((n-1)\mathbb{E}\Bigl [(1-W)^{n-2}\Bigr ] + \theta \Bigr ).
\]
The dual equation is similar to (\ref{dual:5}). Let
\begin{equation*}
h_{\bm{n}}(\bm{x})=\frac{g_{\bm{n}}(\bm{x})}{\mathbb{E}[g_{\bm{n}}(\bm{X})]}
\end{equation*}
then
\begin{equation}
\mathbb{E}_{\bm{X}(0)=\bm{x}}\Bigl [h_{\bm{n}}(\bm{X}(t)\Bigr ] = \mathbb{E}_{\bm{N}(0)=\bm{n}}\Bigl [h_{\bm{N}(t)}(\bm{x})\Bigr ].
\label{dual:5a}
\end{equation}
The multitype death process has transition probabilities which are easy to describe from the sum of the entries $|\bm{N}(t)|$ and (\ref{expansion:0}),
\[
P(\bm{N}(t)=\bm{m} ~\big |~    \bm{N}(0)=\bm{n})
=
\frac{\prod_{j=1}^d{n_j\choose m_j}}{{n\choose m}}
P(|\bm{N}(t)|=m  ~\big |~|\bm{N}(0)|=n).
\]
An equation analogous to the  
 $k$-dimensional Ewens' sampling formula in the Poisson Dirichlet process is to let $\theta_i=\theta/d$, $i=1,\ldots,d$, then the (labelled) sampling formula is
\[
\lim_{d\to \infty}d_{[k]}{n\choose \bm{n}}\mathbb{E}\big [g_{\bm{n}}(\bm{X})\big ],
\]
where $\bm{n}=(n_1,\ldots ,n_k,0,\ldots ,0)$.
The sampling formula limit is
\begin{equation}
\frac{n!\theta^k}{n_1\cdots n_k}
\cdot \frac{\prod_{i=1}^{k}\Bigl [\prod_{j=2}^{n_i}\mathbb{E}\big [(1-W)^{j-2}\big ]\Bigr ]}
{
\prod_{j=1}^{n}\Big [ (j-1)\mathbb{E}\big [(1-W)^{j-2}\big ]+\theta\Big ]
}.
\label{EwensLambda}
\end{equation}
\citet{M2006,L2010} study recursive equations leading to the $\Lambda$-coalescent sampling formula. 
The familiar Ewens' sampling formula is obtained by taking $W=0$.
\section{Acknowledgments}
Robert Griffiths was supported by the Department of Statistics, Stanford University in 2011;
the Institute for Mathematical Sciences at the National University of Singapore in 2011 attending the program on Probability and Discrete Mathematics in Mathematical Biology; the Department of Statistics, University of California, Berkeley, as a Miller Institute Visiting Research Professor in 2012; the Clay Mathematics Institute in a visit to the University of Montreal in 2013;
Banff International Research Station at a meeting \emph{Random Measures and Measure-Valued Processes} in 2013 and the Institute of Statistical Mathematics, Tachikawa, in 2014.
The author thanks the institutions for their support and hospitality. Thanks also to many colleagues and a referee for their helpful discussion, corrections, and suggestions on this research. 
%%%

%%%%%%%%%%%%%%%%%%%%%%%%%%%%%%%

\end{document}